\documentclass{scrartcl}
\usepackage{graphicx}
\usepackage{bm}
\usepackage{amsmath}
\usepackage{amssymb}

\author{Chang, Yifan}
\title{IDENTIFIABILITY OF ELECTRICAL AND HEAT TRANSFER PARAMETERS USING COUPLED BOUNDARY MEASUREMENTS}
\newtheorem{thms}{Theorem}

\newtheorem{pf}{Proof}

\newtheorem{lem}{Lemma}

\begin{document}
\bibliographystyle{plain}
\maketitle

\section{Abstract}
In this paper, we show that a hybrid method using coupled boundary measurements can determine electrical conductivity, thermal conductivity, and the product of heat capacity and heat density within a bounded domain on the plane uniquely up to a boundary-fixing diffromorphism.

\section{Introduction}
\label{intro}
It is a classical problem to try to determine the internal structure of an object by collecting external information. There are many examples of successful approaches to this type of question: CT in medical imaging, seismic imaging in geophysical prospection, and nondestructive testing in industry. In recent years, particularly in medical imaging, there has been interest in developing hybrid methods, i.e. combining different imaging techniques. This idea has proven to be useful, since a given method can provide extra interior information for another which can be used to get a better resolution: see \cite{hyb1}, \cite{hyb2}, \cite{hyb3}, \cite{hyb4}, \cite{hyb5}, \cite{hyb6}.

In \cite{hybrid}, the authors proved uniqueness for a new hybrid method they proposed to determine both the electrical and heat transfer parameters within some bounded domain in $\mathbf{R}^n, n\geq3$, by using coupled boundary measurements. In this paper, we generalize their results to the two dimensional case.

Given a bounded domain $\Omega\subset\mathbf{R}^2$ with smooth boundary, we consider two physical processes. First, we put a voltage distribution $f(\bm{x})\in H^{\frac{1}{2}}(\partial\Omega)$ at the boundary $\partial\Omega$, and assuming there are no sinks or sources
inside $\Omega$, the resulting voltage distribution $u(\bm{x})$ throughout the whole domain is governed by the conductivity equation
\begin{equation}
\label{ce}
\nabla\cdot(\bm{\gamma}(\bm{x})\nabla u(\bm{x}))=0, \quad \bm{x}\in\Omega\qquad\text{and}\quad u|_{\partial\Omega}=f.
\end{equation}
Here $\bm{\gamma}(\bm{x})$ is the electrical conductivity of $\Omega$, the first unknown we want to determine. Since $\bm{\gamma}(\bm{x})$ may be anisotropic, it is represented by a positive definite matrix at each point, and we assume
it is bounded below by a positive number. The standard existence theory assures us that for each $f(\bm{x})\in H^{\frac{1}{2}}(\partial\Omega)$, there exists a unique solution $u(\bm{x})\in H^1(\Omega)$ for (\ref{ce}). We can also let the boundary voltage distribution vary in time, i.e. $f = f(\bm{x},t)$ depends on $t$. For any fixed time $t_0$, we can still solve (\ref{ce}) in the same way to get a voltage distribution $u(\bm{x},t)$ that depends on $t$. 

The second process is heat transfer. It is known that current generates heat which will diffuse throughout the domain. We model this by considering the temperature distribution $\psi(\bm{x})$, which is governed by the heat equations
\begin{eqnarray}
\label{he}
\kappa^{-1}(\bm{x})\partial_t \psi(\bm{x},t)=\nabla\cdot(\bm{A}(\bm{x})\nabla \psi(\bm{x},t))+S(\bm{x},t), \quad \bm{x}\in\Omega\quad\text{and}\quad t\geq 0,\\
\psi(\bm{x},0)=0,\quad \forall \bm{x}\in\Omega\quad \psi(\bm{x},t)=0,\quad \forall \bm{x}\in\partial\Omega\quad\text{and}\quad t\geq 0,
\label{he_bc}
\end{eqnarray}
where $\kappa(\bm{x})=c(\bm{x})^{-1}\rho(\bm{x})^{-1}$ is the reciprocal of the product of the heat capacity $c(\bm{x})$ and density $\rho(\bm{x})$ and $\bm{A}(\bm{x})=(a_{ij}(\bm{x}))$ is the thermal conductivity. The term $S(\bm{x},t)$, which is called the energy density of the electrical field, is given by $S(\bm{x},t)=\nabla u(\bm{x},t)\cdot \bm{\gamma}(\bm{x})\nabla u(\bm{x},t)$, where $u(\bm{x},t)$ is the solution for (\ref{ce}) and will act as a source term in the heat equation. The initial and boundary condition (\ref{he_bc}) means that the temperature is zero (after some shift) at first and the boundary temperature is kept at 0 for all $t\geq 0$. We will assume that the heat transfer is sufficiently slow so that the quasistatic model (\ref{ce}) for the voltage $u(\bm{x})$ is still realistic, as in \cite{hybrid}. 

Define coupled boundary measurements, namely the voltage-to-heat flow map $\Sigma_{\bm{\gamma},\kappa,\bm{A}}$, as follows:
\begin{equation}
\Sigma_{\bm{\gamma},\kappa,\bm{A}}:\quad f(\bm{x},t)\quad\Rightarrow\quad \bm{\nu}\cdot\bm{A}\nabla\psi(\bm{x},t),\quad \bm{x}\in\partial\Omega.
\end{equation}
That is, we set a time-dependent voltage $f(\bm{x},t)$ at the boundary, and measure the out-coming heat flow $\bm{\nu}\cdot\bm{A}\nabla\psi(\bm{x},t)$. We study what information about the internal parameters $\bm{\gamma}$, $\kappa$, and $\bm{A}$, can we recover from the boundary measurements $\Sigma_{\bm{\gamma},\kappa,\bm{A}}$.

The question to determine the inside conductivity by using voltage-to-current measurements at the boundary is known as Calder\'on's inverse problem. It was first proposed by Alberto Calder\'{o}n who came across it while working as an engineer in the 1940's and published his result \cite{c0} in 1980. Since then, Calder\'{o}n's problem has both been applied in industry and become of further theoretical interest. In geophysical prospection, the Schlumberger-Doll company was founded to find oil by using electromagnetic methods. In medical imaging, Calder\'on's inverse problem is known as Electrical Impedance Tomography, which has been used for detecting pulmonary emboli (cf. \cite{eit}). Mathematically, there have been many results on uniqueness (\cite{c0}, \cite{kv0}, \cite{su0}, \cite{sucpam2}, \cite{ht0}, \cite{boazcmp}, \cite{cr}, \cite{nach}, \cite{astala}), stability (\cite{ale}), reconstruction and the corresponding numerical methods (\cite{nachrec}, \cite{mks}, \cite{smi}, \cite{balcpam}). See \cite{review}, \cite{ureview2} for general reviews. For the question of determining the heat parameters through boundary measurements, there have also been some results \cite{heat1}, \cite{heat2}, \cite{heat3}.

Unlike \cite{hybrid}, where the authors consider isotropic conductivity and show that $\Sigma_{\bm{\gamma},\kappa,\bm{A}}$ determines the parameters $\bm{\gamma}$, $\kappa$ and $\bm{A}$ uniquely, we allow the conductivity to be anisotropic, and may encounter some nonuniqueness as a consequence. As we know for the anisotropic Calder\'on's problem, given any boundary-fixing diffeomorphism $\bm{F}$, we can define the pushforward of $\bm{\gamma}$ as $\bm{F}_{*}\bm{\gamma}(\bm{x})=\frac{D\bm{F} \bm{\gamma} D\bm{F}^T}{|D\bm{F}|}\circ \bm{F}^{-1}(\bm{x})$, and this satisfies $\Lambda_{\bm{\gamma}}=\Lambda_{\bm{F}_{*}\bm{\gamma}}$ since the conductivity equation is independent of the choice of coordinates. In the two dimensional case, it has been proved that this will be the only obstacle to uniqueness: for more details, see \cite{anicon}, \cite{lu}, \cite{anisun}, \cite{astalaani}. Similarly we can do the same change of coordinates for the heat equation, so for this new hybrid method, the best one can hope for is that $\Sigma_{\bm{\gamma},\kappa,\bm{A}}$ determines the parameters uniquely up to a boundary-fixing diffeomorphism. We prove that this is indeed the case:
\begin{thms}
\label{thm1}
Assume $\Omega\subset\mathbf{R}^2$ is a bounded domain with smooth boundary, $\kappa_i(\bm{x})\in C^\infty(\bar{\Omega})$, $\kappa_i>c>0$ for some $c$, $\bm{\gamma}_i(\bm{x})$, $\bm{A}_i(\bm{x})\in C^\infty(\bar{\Omega})$ is positive definite and bounded blow by a positive number, $i=1,2$. 
If $\Sigma_{\bm{\gamma}_1,\kappa_1,\bm{A}_1}=\Sigma_{\bm{\gamma}_2,\kappa_2,\bm{A}_2}$, then there exists a smooth diffeomorphism $\bm{F}:\bar{\Omega}\Rightarrow\bar{\Omega}$, which fixes the boundary, i.e. $\bm{F}|_{\partial\Omega}=id$, such that $\kappa_2(\bm{x})=\kappa_1 |D\bm{F}|\circ \bm{F}^{-1}(\bm{x})$, $\bm{\gamma}_2(\bm{x})=\frac{D\bm{F} \bm{\gamma}_1 D\bm{F}^T}{|D\bm{F}|}\circ \bm{F}^{-1}(\bm{x})$ and $\bm{A}_2(\bm{x})=\frac{D\bm{F} \bm{A}_1 D\bm{F}^T}{|D\bm{F}|}\circ \bm{F}^{-1}(\bm{x})$, where $D\bm{F}=(\frac{\partial F_i}{\partial x_j})_{ij}$ is the Jacobian matrix of $\bm{F}$.
\end{thms}

We remark that the regularity assumption is not optimal and may be improved by a more refined analysis. The rest of the paper is organized as follows. In Section \ref{outline} we explain the outline of the proof. In Section \ref{coc} we do a brief review of the nonuniqueness caused by the change of coordinates. In Section \ref{rc}, we show that the conductivity is determined uniquely up to a boundary-fixing diffeomorphism. An important density argument is proved in Section \ref{dense} and after that in Section \ref{rh}, we show how to determine the heat parameters up to the same diffeomorphism as in Section \ref{rc}. From now on, we will use the same regularity assumptions without further mentioning.  

\section{Outline of the Proof}
\label{outline}
In this section, we give the outline of the proof for Theorem \ref{thm1}.
\begin{itemize}
\item First, we recover the quadratic form related to the conductivity equation (\ref{ce}) from the voltage-to-heat flow map $\Sigma_{\bm{\gamma}_i,\kappa_i,\bm{A}_i}$. This is done by taking special input data (almost static ones), and the physical interpretion is that in the static case, the quadratic form stands for the energy needed to maintain the boundary voltage due to the energy lost by the heat effect caused by the current, so this energy should be equal with the heat energy coming out. Then based on the exiting result for the two dimensional Calder\'on's problem, we prove that the two conductivities are the same up to certain boundary-fixing diffeomorphism.
\item The conductivity equation (\ref{ce}) and the heat equation (\ref{he}) are related through the electrical energy density $S_i(\bm{x},t)=\nabla u_i(\bm{x},t)\cdot \bm{\gamma}_i(\bm{x})\nabla u_i(\bm{x},t)$, which is the outcome of (\ref{ce}) and acts as an input of (\ref{he}). Now we consider particularly the case that the boundary voltage $f$ is separated in $\bm{x}$ and $t$, i.e. $f_i(\bm{x},t)=h(\bm{x})g(t)$, then the corresponding solutions will have the form $u_i(\bm{x},t)=u_{i0}(\bm{x})g(t)$, where $u_{i0}$ is the static solution with boundary value $h(\bm{x})$, and $S_i(\bm{x},t)=(\nabla u_{i0}(\bm{x})\cdot \bm{\gamma}_i(\bm{x})\nabla u_{i0}(\bm{x}))g(t)$. Define a bilinear form $B_i(u_i,v_i)\triangleq\nabla u_i(\bm{x})\cdot \bm{\gamma}_i(\bm{x})\nabla v_i(\bm{x})$ where $u_i$, $v_i$ are solutions of (\ref{ce}) with conductivity $\bm{\gamma}_i$. Since (\ref{he}) together with (\ref{he_bc}) are a linear system, having the quadratic form $S_i(\bm{x},t)=B_i(u_{i0},u_{i0})g(t)$ as input is equivalent with having the bilinear form $B(u_{i0},v_{i0})g(t)$ as input. Next we show the space spanned by those bilinear form $B(u_{i0},v_{i0})$ are actually dense in $L^2(\Omega)$(then space spanned by the quadratic form is also dense, we just use the bilinear form for simplicity), which means that $\Sigma_{\bm{\gamma}_i,\kappa_i,\bm{A}_i}$ actually contain information for putting arbitrary variable-separated source $S(\bm{x},t)=w(\bm{x})g(t)$ into ($\ref{he}$). The method to prove density is to use special solutions to the related Schr\"odinger equation(including the Complex Geometric Optic solutions, or CGO solutions for short, and the solutions used by Bukhgeim in \cite{buk}).
\item Based on the fact in the second step and a change of coordinates(resulting from the first step), the problem have become that if for any variable-separated source $S(\bm{x},t)=w(\bm{x})g(t)$, the heat flow coming out are the same, then the heat parameters must coincide. This is proved in \cite{hybrid}. The basic idea is to take special source term(approximately $w(\bm{x})\delta(t)$) and solve the heat equations using eigenfunctions of the operator $\kappa_i \nabla\cdot(\bm{A_i}(\bm{x})\nabla)$ in some weighted $L^2$ space. And with the help of some boundary determination results, it can be proved that all the eigenvalues and eigenfunctions of the two operator $\kappa_i\nabla\cdot(\bm{A_i}(\bm{x})\nabla)$ are the same, which leads to $\bm{A}_1=\bm{A}_2$ then $\kappa_1=\kappa_2$.
\end{itemize}

\section{Obstacle to uniqueness}
\label{coc}
In this section, we review some facts about the obstacle to uniqueness when we allow anisotropic conductivities, i.e. a boundary-fixing diffeomorphism of $\Omega$, or in other words, a change of coordinates and the results in this section are valid for any dimension.

In this whole section, we will assume that $\bm{F}:\bar{\Omega}\Rightarrow\bar{\Omega}$ is a smooth diffeomorphism of $\bar{\Omega}$ which fix the boundary, i.e. $\bm{F}|_{\partial\Omega}=id$. We first start with the conductivity equation (\ref{ce}) and have the following result.

\begin{lem}
\label{coc_ce}
If $u(\bm{x})$ is a solution to (\ref{ce}), then $\tilde{u}(\bm{x})\triangleq u\circ \bm{F}^{-1}(\bm{x})=u(\bm{F}^{-1}(\bm{x}))$ solves (\ref{ce}) with conductivity $\tilde{\bm{\gamma}}=\bm{F}_*\bm{\gamma}\triangleq\frac{D\bm{F} \bm{\gamma} D\bm{F}^T}{|D\bm{F}|}\circ \bm{F}^{-1}$, which is sometimes called the pushforward of $\bm{\gamma}$ and $D\bm{F}=(\frac{\partial F_i}{\partial x_j})_{ij}$ is the Jacobian matrix of $\bm{F}$ as usual. As a result, we have $\Lambda_{\bm{\gamma}}=\Lambda_{\bm{F}_*\bm{\gamma}}$, where $\Lambda_{\bm{\gamma}}$ is the standard Dirichlet-to-Neumann map (DN map) as in the Calder\'on's problem. 
\end{lem}

This can be shown in various ways instead of calculating directly, for example relate the equation to the Laplace-Beltrami operator $\Delta_g$ for some properly chosen metric, i.e. $g\triangleq |\bm{\gamma}|^{\frac{1}{n-2}}\bm{\gamma}^{-1}$, then $\nabla\cdot(\bm{\gamma}(\bm{x})\nabla u(\bm{x}))=|g|^{\frac{1}{2}}\Delta_g u$ and $\tilde{u}(\bm{x})$ solves (\ref{ce}) with conductivity $\tilde{\bm{\gamma}}$ is just another expression in different coordinates for the fact that $u$ solves the Laplace-Beltrami equation which is defined intrinsically. But we may have a problem when $n=2$, since the metric will not be well defined then, we will fix this problem later in the proof of Lemma \ref{coc_he}. Another way to see it is using the fact that $u(\bm{x})$ solves (\ref{ce}) if and only if $u(\bm{x})$ minimize the energy 
$E_{\bm{\gamma}}(u)=\int_\Omega \nabla u(\bm{x})\cdot \bm{\gamma}(\bm{x})\nabla u(\bm{x})\,dx$ within the class of functions with boundary value $f$. And the standard existence and uniqueness results for the solution of (\ref{ce}) implies the existence of a global minimizer. An easy observation using a change of coordinates shows that $E_{\bm{\gamma}}(u(\bm{x}))=E_{\bm{F}_{*}\bm{\gamma}}(u\circ \bm{F}^{-1}(\bm{x}))$. So $u(\bm{x})$ minimize $E_{\bm{\gamma}}(\cdot)$ implies that $u\circ \bm{F}^{-1}$ minimizes $E_{\bm{F}_{*}\bm{\gamma}}(\cdot)$, which means $u\circ \bm{F}^{-1}$ solves (\ref{ce}) with conductivity $\bm{F}_{*}\bm{\gamma}$. To show that $\Lambda_{\bm{\gamma}}=\Lambda_{\bm{F}_*\bm{\gamma}}$, either use the Riemannian geometry approach which will be explained below, or just calculate directly using the explicit expression plus the fact that $\bm{F}$ is identity when restricted to the boundary.

Then we move on to the heat equation and show the following result.

\begin{lem}
\label{coc_he}
If $\psi(\bm{x},t)$ solves the heat equation (\ref{he}) with the initial and boundary value condition (\ref{he_bc}). Define $\tilde{\kappa}(\bm{x})\triangleq |D\bm{F}(\bm{F}^{-1}(\bm{x}))|\kappa(\bm{F}^{-1}(\bm{x}))$, $\tilde{\psi}(\bm{x},t)\triangleq \psi(\bm{F}^{-1}(\bm{x}),t)$, $\tilde{S}(\bm{x},t)\triangleq |D\bm{F}(\bm{F}^{-1}(\bm{x}))|^{-1}S(\bm{F}^{-1}(\bm{x}),t)$, $\tilde{A}(\bm{x})=\bm{F}_*\bm{A}(\bm{x})\triangleq\frac{D\bm{F} \bm{A} D\bm{F}^T}{|D\bm{F}|}\circ \bm{F}^{-1}(\bm{x})$, then $\tilde{\psi}(\bm{x},t)$ solves the system with parameters $\tilde{\kappa}(\bm{x})$, $\tilde{\bm{A}}(\bm{x})$ and $\tilde{S}(\bm{x},t)$. As a result, the two system give the same out-coming heat flow at the boundary, i.e. $\bm{\nu}\cdot\bm{A}\nabla\psi(\bm{x},t)=\bm{\nu}\cdot\tilde{\bm{A}}\nabla\tilde{\psi}(\bm{x},t)$, $\bm{x}\in\partial\Omega$. 
\end{lem}

\begin{pf}
Denote that $\bm{y}=\bm{F}(\bm{x})$, according to the definition, the only thing need to show here is that
\begin{equation}
|D\bm{F}(\bm{x})|^{-1}\nabla_{\bm{x}}\cdot(\bm{A}(\bm{x})\nabla_{\bm{x}} u(\bm{x}))=\nabla_{\bm{y}}\cdot(\tilde{\bm{A}}({\bm{y}})\nabla_{\bm{y}} \tilde{u}(\bm{y}))=
\nabla_{\bm{y}}\cdot(\bm{F}_*\bm{A}({\bm{y}})\nabla_{\bm{y}} u(\bm{F}^{-1}(\bm{y})))
\end{equation}
We will use the metric arguments mentioned above. First when $n\geq 3$, $\forall \bm{A}(\bm{x})$ positive definite with a positive lower bound, we associate it with a metric 
\begin{equation}
\label{metric}
g_{\bm{A}}(\bm{x})\triangleq |\bm{A}(\bm{x})|^{\frac{1}{n-2}}\bm{A}(\bm{x})^{-1}
\end{equation}
And here are some basic facts which can be verified easily.
\begin{itemize}
\item $|g_{\bm{A}}|=|\bm{A}|^{\frac{2}{n-2}}$
\item $\nabla\cdot(\bm{A}(\bm{x})\nabla u(\bm{x}))=|g_{\bm{A}}|^{\frac{1}{2}}\Delta_{g_{\bm{A}}}u=|\bm{A}|^{\frac{1}{n-2}}\Delta_{g_{\bm{A}}}u$
\item The metric associated with $\tilde{\bm{A}}=\bm{F}_{*}\bm{A}$ is just the pullback of $g_{\bm{A}}$ under $\bm{F}^{-1}$, i.e. $g_{\bm{F}_{*}\bm{A}}=(\bm{F}^{-1})^{*}g_{\bm{A}}$. $|\bm{F}_{*}\bm{A}(\bm{y})|=|D\bm{F}(\bm{x})|^{2-n}|\bm{A}(\bm{x})|$, $|g_{\bm{F}_*\bm{A}}(\bm{y})|=|D\bm{F}(\bm{x})|^{-2}|g_{\bm{A}}(\bm{x})|$
\end{itemize}
Then it is easy to see that
\begin{equation*}
\begin{split}
\nabla_{\bm{x}}\cdot(\bm{A}(\bm{x})\nabla_{\bm{x}} u(\bm{x}))
&=|g_{\bm{A}}|^{\frac{1}{2}}\Delta_{g_{\bm{A}}}u
=(\frac{|g_{\bm{A}}|}{|(\bm{F}^{-1})^{*}g_{\bm{A}}|})^{\frac{1}{2}}
|(\bm{F}^{-1})^{*}g_{\bm{A}}|^{\frac{1}{2}}
\Delta_{(\bm{F}^{-1})^{*}g_{\bm{A}}}u\circ\bm{F}^{-1}\\
&=(\frac{|g_{\bm{A}}|}{|g_{\tilde{\bm{A}}}|})^{\frac{1}{2}}
|g_{\tilde{\bm{A}}}|^{\frac{1}{2}}\Delta_{g_{\tilde{\bm{A}}}}u\circ\bm{F}^{-1}
=|D\bm{F}(\bm{x})|\nabla_{\bm{y}}\cdot(\tilde{\bm{A}}(\bm{y})\nabla_{\bm{y}} \tilde{u}(\bm{y}))
\end{split}
\end{equation*}

To show that $\bm{\nu}\cdot\bm{A}\nabla\psi(\bm{x},t)=\bm{\nu}\cdot\tilde{\bm{A}}\nabla\tilde{\psi}(\bm{x},t)$, we can either calculate directly, or relate it with the interior product of $\nabla_g u$ and the volume form $d_gV$ since we have the relation
\begin{equation}
\label{volume}
(\nabla_g \psi\,\lrcorner\, d_g V)|_{\partial\Omega}=(\bm{\nu}\cdot\bm{A}\nabla\psi(\bm{x},t))dS,
\end{equation}
where $dS$ is the Euclidean volume form restricted to $\partial\Omega$. Then we get what we want since $\nabla_g \psi\,\lrcorner\, d_g V$ is independent of the choice of coordinates plus the fact that $\bm{F}$ is identity at the boundary.

When $n=2$, we can not use the arguments above since the term $|\bm{A}(\bm{x})|^{\frac{1}{n-2}}$ in (\ref{metric}) makes no sense, or from another perspective, in the Laplace-Beltrami operator 
$\Delta_g=|g|^{-\frac{1}{2}}\nabla\cdot(|g|^{\frac{1}{2}}g^{-1}\nabla)$, the matrix $|g|^{\frac{1}{2}}g^{-1}$ always has determinant one for the two dimensional case, which is not satisfied by general $\bm{A}$. We fix this by normalize $\bm{A}$ first, define $\bm{A}(\bm{x})=|\bm{A}(\bm{x})|^{\frac{1}{2}}(\frac{1}{|\bm{A}(\bm{x})|^{\frac{1}{2}}}\bm{A}(\bm{x}))\triangleq|\bm{A}(\bm{x})|^{\frac{1}{2}}\bm{A}_0(\bm{x})$, so $|\bm{A}_0|=1$ and 
$$\nabla\cdot(\bm{A}\nabla u)=\nabla\cdot(|\bm{A}|^{\frac{1}{2}}\bm{A}_0\nabla u)
=|\bm{A}|^{\frac{1}{2}}\nabla\cdot(\bm{A}_0\nabla u)+\nabla |\bm{A}|^{\frac{1}{2}}\cdot(\bm{A}_0\nabla u).$$ Then just simply take a metric $g=\bm{A}^{-1}$, we get
$$\nabla\cdot(\bm{A}(\bm{x})\nabla u(\bm{x}))=\Delta_g u+<\nabla_g \log{|\bm{A}|^{\frac{1}{2}}},\nabla_g u>_g,$$
and similarly, $\nabla\cdot(\tilde{\bm{A}}(\bm{y})\nabla \tilde{u}(\bm{y}))=\Delta_{\tilde{g}} \tilde{u}+<\nabla_{\tilde{g}} \log{|\tilde{\bm{A}}|^{\frac{1}{2}}},\nabla_{\tilde{g}} \tilde{u}>_{\tilde{g}}$. Notice that when $n=2$, $\tilde{g}\triangleq\tilde{\bm{A}}^{-1}$ does not coincides with $(\bm{F}^{-1})^*(g)$, but satisfies that $\tilde{g}(\bm{y})=|D\bm{F}(\bm{x})|(\bm{F}^{-1})^* g(\bm{y})$. Combined the fact that $|\bm{A}(\bm{x})|=|\tilde{\bm{A}}(\bm{y})|$, we have 
$$\Delta_{\tilde{g}} \tilde{u}(\bm{y})=|D\bm{F}(\bm{x})|^{-1}\Delta_g u(\bm{x}),\, <\nabla \log{|\tilde{\bm{A}}|^{\frac{1}{2}}},\nabla \tilde{u}>_{\tilde{g}}|_{\bm{y}}=|D\bm{F}|^{-1}<\nabla \log{|\bm{A}|^{\frac{1}{2}}},\nabla u>_g|_{\bm{x}},$$ which finishes the proof that $\tilde{\psi}$ solves the equation with $\tilde{\kappa}$, $\tilde{\bm{A}}$ and $\tilde{S}$. 

To show $\bm{\nu}\cdot\bm{A}\nabla\psi(\bm{x},t)=\bm{\nu}\cdot\tilde{\bm{A}}\nabla\tilde{\psi}(\bm{x},t)$, we can also either calculate directly or use the same arguments mentioned above with some small modification. Notice that when $n=2$, the relation (\ref{volume}) becomes
\begin{equation}
|g|^{-\frac{1}{2}}(\nabla_g \psi\,\lrcorner\, d_g V)|_{\partial\Omega}=(\bm{\nu}\cdot\bm{A}\nabla\psi(\bm{x},t))dS.
\end{equation}
Then use the fact $\tilde{g}(\bm{y})=|D\bm{F}(\bm{x})|(\bm{F}^{-1})^* g(\bm{y})$, $|\tilde{g}(\bm{y})|=|g(\bm{x})|$, $d_{(\bm{F}^{-1})^* g} V=|D\bm{F}|^{-1}d_{\tilde{g}} V$, $\nabla_{(\bm{F}^{-1})^* g} \tilde{\psi}=|D\bm{F}|\nabla_{\tilde{g}}\tilde{\psi}$, we get 
$$\nabla_g \psi\,\lrcorner\, d_g V|_{\bm{x}}=(\nabla_{(\bm{F}^{-1})^* g} \psi\,\lrcorner\, d_{(\bm{F}^{-1})^* g} V)|_{\bm{y}}=\nabla_{\tilde{g}} \psi\,\lrcorner\, d_{\tilde{g}} V|_{\bm{y}},$$
plus the fact that $\bm{F}$ is identity at the boundary, we finish the proof.
\end{pf}

\section{Determination of the Conductivity}
\label{rc}
In this section, we recover the conductivity from $\Sigma_{\bm{\gamma},\kappa,\bm{A}}$, most of the arguments follow from \cite{hybrid}. We start with a property that will be used to solve the heat equation. Consider the operator $P\triangleq -\kappa\nabla\cdot(\bm{A}\nabla)$ in the weighted $L^2$ space $L^2_{\kappa^{-1}}(\Omega)$ with weight function $\kappa^{-1}(x)$, i.e. the space equipped with the weighted inner product $<u,v>_{\kappa^{-1}}=\int_\Omega u\bar{v}\kappa^{-1}\,d\bm{x}$, we start with the domain of $P$ as $C_0^\infty(\Omega)$ and then extend it to $H^1_0(\Omega)\cap H^2(\Omega)$. It is well known that, given $\kappa$, $\bm{A}$ both positive with lower bounds strictly greater than zero, $P$ is positive self-adjoint, furthermore the spectrum of $P$ consist of real positive eigenvalues $\{\lambda_i\}_{i=1}^\infty$ of finite multiplicity which accumulating at $+\infty$(we may assume that $\lambda_i$ is non-decreasing), and the corresponding eigenfunctions $\{\phi_i\}_{i=1}^\infty$ form an orthonormal basis of $L^2_{\kappa^{-1}}(\Omega)$.

And now we return to the topic of this section, determining the conductivity through $\Sigma_{\bm{\gamma},\kappa,\bm{A}}$. The method used here have a strong physical interpretation. Recall that according to the classical theory of Calder\'on problem, the conductivity is determined by the Dirichlet-to-Neumann map $\Lambda_{\bm{\gamma}}$ up to a boundary-fixing diffeomorphism. And knowing the Dirichlet-to-Neumann map $\Lambda_{\bm{\gamma}}$ is actually equivalent with knowing the quadratic form, or the energy function  $Q_{\bm{\gamma}}(f)=E_{\bm{\gamma}}(u)=\Lambda_{\bm{\gamma}}f(f)=\int_{\partial\Omega}\bm{\gamma}\nabla u\cdot \bm{n}\, dS=\int_\Omega \nabla u\cdot\bm{\gamma}\nabla u\, d\bm{x}$, where $u$ is the solution to (\ref{ce}) with boundary value $f$. The physical meaning of $Q(f)=E_{\bm{\gamma}}(u)$ is that, it is the energy needed to maintain the boundary voltage $f$. We know that energy is conservative, so it must transfer into some other form, and the answer is heat. So it is quite nature to expect that in the static case, i.e. when we fix the boundary voltage to be $f$ and wait long enough so that the temperature does not change, the energy we put into the system should equal to the heat coming out, i.e. $Q(f)=E_{\bm{\gamma}}(u)=-\int_{\partial\Omega}\bm{A}\nabla \psi \cdot \bm{n}\,dS=-\int_{\partial\Omega} \Sigma_{\bm{\gamma},\kappa,\bm{A}}f \,dS$. This can be proved as follows, assume $u$ and $\psi$ solves the static conductivity and heat equation, i.e. $u$, $v$ solve (\ref{ce}), (\ref{he}) and don't depend on $t$, then $\int_{\partial\Omega}\bm{A}\nabla \psi \cdot \bm{n}\,dS=\int_{\Omega}\nabla\cdot(\bm{A}\nabla \psi)\,d\bm{x}=-\int_{\Omega} S(\bm{x})\,d\bm{x}=-\int_{\Omega}\nabla u \cdot\bm{\gamma}\nabla u \,d\bm{x}$. This is just an intuitive proof since in this case $\psi$ doesn't satisfy the initial condition. But it gives us the essential idea and now we prove it rigorously. The idea is still the same, we try to recover the quadratic form $Q_{\bm{\gamma}}$ by taking the static input data $f(\bm{x},t)=f(\bm{x})$, and show that $\int_{\partial\Omega}(\Sigma_{\bm{\gamma},\kappa,\bm{A}}f)(\bm{x},t)\, dS\rightarrow -Q_{\bm{\gamma}}(f)$, as $t\rightarrow +\infty$ 
\begin{lem}
If we take $f(\bm{x},t)=f(\bm{x})$, then $\lim_{t\rightarrow +\infty}\int_{\partial\Omega}(\Sigma_{\bm{\gamma},\kappa,\bm{A}}f)(\bm{x},t)\, dS=-Q_{\bm{\gamma}}(f)$.
\end{lem}
\begin{pf}
Assume $u_0(\bm{x})$ solves (\ref{ce}) with the static boundary data $f(\bm{x},t)=f(\bm{x})$, then the source term will also be static (independent of $t$), $S(\bm{x},t)=S(\bm{x})=\nabla u_0\cdot\bm{\gamma}\nabla u_0$ and $Q_{\bm{\gamma}}(f)=\int_\Omega S(\bm{x})\, d\bm{x}$. After that we solve the heat equation (\ref{he}) by the decomposition $\psi(\bm{x},t)=\psi_0(\bm{x})+\psi_1(\bm{x},t)$, where $\psi_0(\bm{x})$ solves the static heat equation with the static source term,
\begin{equation}
\nabla\cdot(\bm{A}\nabla \psi_0(\bm{x}))+S(\bm{x})=0,\,\bm{x}\in\Omega,\quad \psi_0(\bm{x})=0,\quad \forall \bm{x}\in\partial\Omega,
\label{psi0}
\end{equation}
and $\psi_1(\bm{x},t)$ fix the initial condition,
\begin{equation}
\kappa^{-1}\partial_t \psi_1=\nabla\cdot(\bm{A}\nabla \psi_1),\, \psi_1(\bm{x},0)=-\psi_0(\bm{x}),\, \forall \bm{x}\in\Omega,\quad \psi_1(\bm{x},t)=0,\,\forall \bm{x}\in\partial\Omega.
\label{psi1}
\end{equation}
From the standard existence theory, we know there is always a unique $\psi_0(\bm{x})\in H_0^1(\Omega)$ solves (\ref{psi0}). Also notice that 
\begin{equation*}
\int_{\partial\Omega}\bm{A}\nabla \psi_0 \cdot \bm{n}\,dS=\int_\Omega \nabla\cdot(\bm{A}\nabla \psi_0) \,d\bm{x}
=-\int_\Omega S(\bm{x})\,d\bm{x}=-Q_{\bm{\gamma}}(f).
\end{equation*}
And for $\phi_1$, we solve it by using the eigenfunction expansion as follows. First we do the eigenfunction decomposition for the initial value $-\psi_0$,
\begin{equation}
-\psi_0(\bm{x})=\sum_{i=1}^\infty a_i\phi_i(\bm{x}),\quad a_i=<-\psi_0,\phi_i>_{\kappa^{-1}}=-\int_\Omega \psi_0\bar{\phi}_i\kappa^{-1}\,d\bm{x}.
\end{equation} 
Plug this into the equation, we get $\psi_1(\bm{x},t)=\sum_{i=1}^\infty a_i e^{-\lambda_i t} \phi_i(\bm{x})$. Since $\lambda_i\geq\lambda_1>0$, we know that $\psi_1$ converge to $0$ exponentially as $t\rightarrow +\infty$, so is the out-coming heat flow for $\psi_1$. Then we have $$\lim_{t\rightarrow +\infty}\int_{\partial\Omega}(\Sigma_{\bm{\gamma},\kappa,\bm{A}}f)(\bm{x},t)\, dS=\int_{\partial\Omega}\bm{A}\nabla \psi_0 \cdot \bm{n}\,dS=-Q_{\bm{\gamma}}(f)$$.  
\end{pf}
In this way we recover the quadratic form $Q_{\bm{\gamma}}$, then according to the theory for the two dimensional Calder\'on's problem (see \cite{astala}, \cite{astalaani} for reference), this determine the conductivity up to a boundary-fixing diffeomorphism, and we also have some recovering algorithms for it.
\begin{lem}
\label{rc2}
If $\Sigma_{\bm{\gamma}_1,\kappa_1,\bm{A}_1}=\Sigma_{\bm{\gamma}_2,\kappa_2,\bm{A}_2}$, then there is a boundary-fixing diffeomorphism $\bm{F}$ of $\bar{\Omega}$, s.t. $\bm{\gamma}_2=\bm{F}_*\bm{\gamma}_1$
\end{lem}

\section{A Density Argument}
\label{dense}
In this section, we prove the density argument mentioned in Section \ref{intro} as a preparation for determining the heat parameters. The main result is stated as follows.
\begin{lem}
\label{den}
Assume $\bm{\gamma}$ is regular enough, then the space spanned by $\nabla u(\bm{x})\cdot \bm{\gamma}(\bm{x})\nabla v(\bm{x})$, where $u(\bm{x})$ and $v(\bm{x})$ are arbitrary solutions of (\ref{ce}), is dense in $L^2(\Omega)$. 
\end{lem}

In order to show $span\{\nabla u(\bm{x})\cdot \bm{\gamma}(\bm{x})\nabla v(\bm{x})\}$ is dense, we prove that $span\{\nabla u(\bm{x})\cdot \bm{\gamma}(\bm{x})\nabla v(\bm{x})\}^{\perp}={0}$, i.e. if $f(\bm{x})\in L^2(\Omega)$ satisfies that $\int_\Omega f(\bm{x})\nabla u(\bm{x})\cdot \bm{\gamma}(\bm{x})\nabla v(\bm{x})\, d\bm{x}=0$ for any pair of solutions $u(\bm{x})$ and $v(\bm{x})$, then $f=0$.

First we claim that since now we are working in the two dimensional case, we may assume that $\bm{\gamma}$ is isotropic, i.e. $\bm{\gamma}=\gamma(\bm{x})I$. Since it have been proved in \cite{astalaani} that $\forall \bm{\gamma}$, there is a quasiconformal homeomorphism $\bm{F}\in W_{loc}^{1,p}(\mathbb{C};\mathbb{C})$ of the whole plane with asymptotic behaviour $\bm{F}(z)=z+\mathcal{O}(\frac{1}{z})$, such that $\bm{F}_*\bm{\gamma}=det(\bm{\gamma}\circ \bm{F}^{-1})^{\frac{1}{2}}I$, i.e. we can find a change of coordinates to convert any anisotropic conductivity to a isotropic one. And it is not hard to see that $f\perp span\{\nabla u(\bm{x})\cdot \bm{\gamma}(\bm{x})\nabla v(\bm{x})\}$ is equivalent with $(f\circ\bm{F}^{-1})\perp span\{\nabla (u\circ\bm{F}^{-1})\cdot \bm{F}_*\bm{\gamma}(\bm{x})\nabla (v\circ\bm{F}^{-1})\}$, so the density of $span\{\nabla u(\bm{x})\cdot \bm{\gamma}(\bm{x})\nabla v(\bm{x})\}$ in $L^2(\Omega)$ is equivalent with the density of $span\{\nabla (u\circ\bm{F}^{-1})\cdot \bm{F}_*\bm{\gamma}(\bm{x})\nabla (v\circ\bm{F}^{-1})\}$ in $L^2(\bm{F}(\Omega))$.

To show that $f=0$, we use some special solutions of (\ref{ce}). The outline is that we first use the CGO solutions to show that $f\chi_\Omega\in H^s(\mathbf{R}^2)$ for any $s>0$. Then by the Sobolev embedding, we know that $f$ must vanish to infinite order at the boundary $\partial\Omega$(actually we only need it to vanish to the second order). Then we can apply integration by parts $0=\int_\Omega f\nabla u\cdot \gamma\nabla v\, d\bm{x}=\frac{1}{2}\int_\Omega f\nabla\cdot (\gamma\nabla(uv))\, d\bm{x}=\frac{1}{2}\int_\Omega \nabla\cdot (\gamma\nabla(f))uv\, d\bm{x}$. After that, use the Liouville transform and the solutions Bukhgeim proposed in \cite{buk}, we derive that $\frac{1}{\gamma}\nabla\cdot (\gamma\nabla(f))=0$ in $\Omega$. Combining the facts that $f$ vanishes to all orders at the boundary (once again, we only need it to vanish to the first order), the unique continuation results for elliptic equations show that $f=0$ which finishes the proof.

So now we focus on the first key ingredient in the proof, prove the smoothness of $f\chi_\Omega$ by using CGO solutions.
\begin{lem}
\label{smooth}
Assume $\gamma$ is smooth enough, then if $f\in L^2(\Omega)$, $f\perp span\{\nabla u(\bm{x})\cdot \gamma(\bm{x})\nabla v(\bm{x})\}$, where $u$, $v$ are arbitrary solutions to (\ref{ce}), then $f\chi_\Omega\in H^s(\mathbf{R}^2)$ for any $s>0$.  
\end{lem}
Recall the standard Liouville transform that transform (\ref{ce}) into a Schr\"odinger equation with potential $q=\frac{\Delta \sqrt{\gamma}}{\sqrt{\gamma}}$. In \cite{sugeneric}, the authors proved that if $g\perp span\{u v\}$, where $u$ and $v$ are the solutions to the corresponding Schr\"odinger equation, then $g\chi_\Omega\in H^s(\mathbf{R}^2)$ for $0<s<1$. They used pairs of CGO solutions whose product was approximately $e^{i\bm{k}\cdot\bm{x}}$ to derive estimates for the Fourier transform of $g\chi_\Omega$ for large $\bm{k}$. We extend this result and show that $g\chi_\Omega$ actually belongs to $H^s(\mathbf{R}^2)$ for any $s>0$ and similar arguments can be used to show that $f\chi_\Omega\in H^s(\mathbf{R}^2)$, $\forall s>0$ when we come back to the conductivity equation.

First we review and modify the construction of CGO solutions for the Schr\"odinger equation, which was first introduced in \cite{su1}. For $\forall\bm{\eta}\in\mathbb{C}^2$ large enough, s.t. $\bm{\eta}\cdot \bm{\eta}=0$, or equivalently $\bm{\eta}=\frac{1}{2}(\pm\bm{k}^\perp +i\bm{k})$ (the factor $\frac{1}{2}$ is just for notation simplicity), $\forall \bm{k}\in\mathbf{R}^2$ with sufficiently large norm and $\bm{k}^\perp$ stands for the vector obtained by rotating $\bm{k}$ clockwise by $\frac{\pi}{2}$, we construct CGO solution $u$ for the Schr\"odinger equation in the form $u(\bm{x},\bm{k})=e^{\bm{\eta}\cdot\bm{x}}(1+r(\bm{x},\bm{k}))$. The correction term $r(\bm{x},\bm{k})$ can be proved to be small in certain norm (like the $W^{2,p}(\Omega)$ norm) by solving $r(\bm{x},\bm{k})$ out explicitly. Actually it was shown in \cite{su1} that for $\bm{\eta}=\frac{1}{2}(\bm{k}^\perp +i\bm{k})$, we have an expansion for the correction term $r(\bm{x},\bm{k})$,
\begin{equation}
\label{expansion}
r(\bm{x},\bm{k})=\frac{a(\bm{x})}{i\bar{k}}+\frac{b(\bm{x},\bm{k})}{(i\bar{k})^2},
\end{equation}
where we have uniform bounds(in $\bm{k}$) for the norm of $b(\bm{x},\bm{k})$. And based on this, \cite{sugeneric} showed that $g\chi_\Omega\in H^s(\mathbf{R}^2)$, $0<s<1$ if $g\perp span\{u_s v_s\}$. To extend the result to any $s>0$, we need to modify the expansion (\ref{expansion}) with the following result.
\begin{lem}
\label{expand}
$\forall n\in\mathbf{N}_+$, when the potential $q$ is regular enough(like in $C^{n+2}(\bar{\Omega})$) and $|\bm{k}|$ is large enough(depends on $||q||_{C^{n+2}(\bar{\Omega})}$), we can construct CGO solutions with correction term $r(\bm{x},\bm{k})$ in the following form,
\begin{equation}
\label{higherorder}
r(\bm{x},\bm{k})=\sum_{j=1}^{n-1} \frac{a_j(\bm{x})}{(i\bar{k})^j}+\frac{b_n(\bm{x},\bm{k})}{(i\bar{k})^n},
\end{equation}
where $a_j=(-2\partial+\chi Pq)^{j-1}\chi Pq$, $\chi(\bm{x})$ is any smooth cut-off function which equals $1$ in $\Omega$, the norm of $b$ (like the $W^{2,p}(\Omega)$ norm) is bounded by the norm of $q$ (like the $W^{2,p}(\Omega)$ norm) and independent of $\bm{k}$.  
\end{lem}

In the following, we will use the usual notation $\partial=\frac{1}{2}(\partial_x-i\partial_y)$, $\bar{\partial}=\frac{1}{2}(\partial_x+i\partial_y)$ and switch between vector and complex notation without further mentioning, for example $\bm{k}$, $\bm{x}$ stands for the vector $(k_1,k_2)^T$, $(x_1,x_2)^T$ while $k$, $z$ stands for the complex number $k_1+ik_2$, $x_1+ix_2$. Then $e^{\bm{\eta}\cdot\bm{x}}(1+r(\bm{x},\bm{k}))=e^{\frac{i}{2}\bar{k}z}(1+r(\bm{x},\bm{k}))$ solves the Schr\"odinger equation $(-\Delta+q)u=0$ is equivalent with $(\Delta+2\bm{\eta}\cdot \nabla-q)r=q$, or written in complex notation $(4\bar{\partial}\partial+2i\bar{k}\bar{\partial}-q)r=q$. This is solved by first formally setting $r=\sum_1^\infty r_i$, where $r_i$ satisfies $(4\bar{\partial}\partial+2i\bar{k}\bar{\partial})r_{i}=qr_{i-1}$($r_0=1$), then proving the summation converges(this is where we need $|\bm{k}|$ to be large enough) and solves the equation. So the core of the problem have reduced to solve the following equation inside the domain,
\begin{equation}
(4\bar{\partial}\partial+2i\bar{k}\bar{\partial})r=f.
\end{equation}
The first thing to do is factorize the equation like $2\bar{\partial}(2\partial+i\bar{k})r=f$. Instead of solving it inside the domain, we extend $f$ and do all the operations on the whole plane while only requiring the constructed solution satisfies the equation inside $\Omega$. Notice that, there are various way to extend $f$, later on we will actually need the extension to be regular enough. Then we can first invert the $\bar{\partial}$ equation by the standard Cauchy transform $P$, i.e. $Pf(z)=\frac{1}{2\pi}f*\frac{1}{\cdot}(z)=\frac{1}{2\pi}\int_{\mathbf{R}^2}\frac{f(w)}{z-w}\, dw$, where $dw$ is the usual Lebesgue measure on the plane (we have modified the standard definition by a factor of $\frac{1}{2}$ for simplicity). The properties of $P$ have been studied quite carefully, see \cite{cauchy} for instance. If $f$ is regular enough, we know $Pf$ is always one more time differentiable than $f$ and then $Pf$ solves the equation in the weak sense implies $Pf$ also solves the equation in classical sense point-wisely. So as mentioned above, since we only require the solution to solve the equation inside the domain, we may multiply $Pf$ by a smooth cut-off function $\chi(\bm{x})$ which is $1$ inside $\Omega$.

Now we are trying to solve the equation 
\begin{equation}
\label{coreequ}
(2\partial+i\bar{k})r=g,
\end{equation}
where we are taking $g=\chi Pf$. Notice that $2\partial(e^{i\bm{k}\cdot\bm{x}})=i\bar{k}e^{i\bm{k}\cdot\bm{x}}$, similar to integration factor, we have $(2\partial+i\bar{k})r=e^{-i\bm{k}\cdot\bm{x}}2\partial e^{i\bm{k}\cdot\bm{x}}r$, which means that $2\partial+i\bar{k}$ is actually $2\partial$ conjugated by $e^{i\bm{k}\cdot\bm{x}}$. Then we can solve the equation (\ref{coreequ}) using $\bar{P}$, which is defined as the convolution with $\frac{1}{\bar{z}}$ and acts as the inverse of $2\partial$, and get 
\begin{equation}
r=e^{-i\bm{k}}\bar{P}(e^{i\bm{k}}g)=e^{-i\bm{k}}\bar{P}(e^{i\bm{k}}\chi Pf).
\end{equation}
Due to the mapping properties of $P$ and $\bar{P}$, if we solve the equation this way, the $L^p$ norm of $r$ is the same order as the size of $f$ and independent of $\bm{k}$.

But the expansion structure is not very clear if we solve it like this. Instead since (\ref{coreequ}) is a linear equation with constant coefficient, we may solve it using Fourier transform and get $i(\bar{k}+\bar{w})\hat{r}(w)=\hat{g}(w)$, then $\hat{r}(w)=\frac{1}{i\bar{k}+i\bar{w}}\hat{g}(w)$. Rewrite 
\begin{equation}
\frac{1}{i\bar{k}+i\bar{w}}=\frac{1}{i\bar{k}}\frac{1}{1+\frac{i\bar{w}}{i\bar{k}}}=\frac{1}{i\bar{k}}(1+(-\frac{i\bar{w}}{i\bar{k}})+...+(-\frac{i\bar{w}}{i\bar{k}})^{n-1}+(-\frac{i\bar{w}}{i\bar{k}})^{n}\frac{1}{1+\frac{i\bar{w}}{i\bar{k}}}),
\end{equation}
then we have the expansion
\begin{equation}
\hat{r}(w)=\frac{\hat{g}(w)}{i\bar{k}}+\frac{-i\bar{w}\hat{g}(w)}{(i\bar{k})^2}+...+\frac{(-i\bar{w})^{n-1}\hat{g}(w)}{(i\bar{k})^{n}}+\frac{1}{(i\bar{k})^{n}}\frac{(-i\bar{w})^n\hat{g}(w)}{i\bar{k}+i\bar{w}}.
\end{equation}
Remember that $\widehat{2\partial g}(w)=i\bar{w}\hat{g}(w)$, so it is clear from here that why we require some regularity for $g=\chi Pf$, which can be obtained if we require some regularity for $f$. Then if we take the inverse Fourier transform, we have
\begin{equation}
r(\bm{x},\bm{k})=\frac{g(\bm{x})}{i\bar{k}}+\frac{-2\partial g(\bm{x})}{(i\bar{k})^2}+...+\frac{(-2\partial)^{n-1}g(\bm{x})}{(i\bar{k})^{n}}+\frac{r^{(n)}(\bm{x},\bm{k})}{(i\bar{k})^{n}},
\end{equation}
where $r^{(n)}(\bm{x},\bm{k})$ satisfies the equation $(2\partial+i\bar{k})r^{(n)}=(-2\partial)^{n}g=(-2\partial)^{n}\chi Pf$. Based on the argument above, the $L^p$ norm of $r^{(n)}$ can be controlled by the norm of $D^{(n-1)}f$ and is independent of the size of $\bm{k}$. To sum up, the more regularity we have in $f$, the higher order we can expand in (\ref{higherorder}). We remark here that of course we can get the same result here by setting 
\begin{equation}
r=\sum_{j=1}^{n}\frac{(-2\partial)^{j-1} g}{(i\bar{k})^j}+\frac{r^{(n)}}{(i\bar{k})^{n}},
\end{equation}
directly and another thing is we will get the same solution no matter which order we expand to, as long as they make sense, i.e. $\forall n, m\in\mathbf{N}_+$,
\begin{equation}
\begin{split}
r(\bm{x},\bm{k})&=\frac{g(\bm{x})}{i\bar{k}}+\frac{-2\partial g(\bm{x})}{(i\bar{k})^2}+...+\frac{(-2\partial)^{n-1}g(\bm{x})}{(i\bar{k})^{n}}+\frac{r^{(n)}(\bm{x},\bm{k})}{(i\bar{k})^{n}}\\
& =\frac{g(\bm{x})}{i\bar{k}}+\frac{-2\partial g(\bm{x})}{(i\bar{k})^2}+...+\frac{(-2\partial)^{m-1}g(\bm{x})}{(i\bar{k})^{m}}+\frac{r^{(m)}(\bm{x},\bm{k})}{(i\bar{k})^{m}},
\end{split}
\end{equation}

Now to get the expansion structure mentioned in Lemma \ref{expand}, we first regularly extend the potential $q$ to the whole plane with compact support by standard technique($\forall m\in\mathbf{N}_+$, we can make $q_{ext}\in C_0^m(\mathbf{R}^2)$ and we will still use $q$ instead of $q_{ext}$ for simplicity), then apply the above arguments to $r_i$ up to $(n-i)$-th order, $0\leq i\leq n-2$ and sum them up. The $i$-th order term has the form $(-2\partial+\chi Pq)^{i-1}\chi Pq$ can be proved by induction. Finally the summation $r=\sum_{i=0}^{\infty}r_i$ converges and solves the equation can be proved by the same arguments used in \cite{su1} which we are not going to repeat here since basically we are talking about the same solutions. For instance, to get an expansion up to the third order, we have
\begin{equation}
\begin{split}
r_0 =&\frac{\chi Pq}{i\bar{k}}+\frac{-2\partial \chi Pq}{(i\bar{k})^2}+\frac{(-2\partial)^{2}\chi Pq}{(i\bar{k})^{3}}+\frac{r_0^{(3)}}{(i\bar{k})^{3}}\\
r_1 =&\frac{\chi P(qr_0)}{i\bar{k}}+\frac{-2\partial \chi P(qr_0)}{(i\bar{k})^2}+\frac{r_1^{(2)}}{(i\bar{k})^{2}}\\
=&\frac{1}{i\bar{k}}\chi P(q\frac{\chi Pq}{i\bar{k}}+q\frac{-2\partial \chi Pq}{(i\bar{k})^2}+q\frac{r_0^{(2)}}{(i\bar{k})^{2}})+\\
&\frac{1}{(i\bar{k})^2}(-2\partial)\chi P(q\frac{\chi Pq}{i\bar{k}}+q\frac{r_0^{(1)}}{i\bar{k}})+\frac{r_1^{(2)}}{(i\bar{k})^2},\\
\end{split}
\end{equation}
where $r_0^{(i)}$ solves $(2\partial+i\bar{k})r_0^{(i)}=(-2\partial)^{i}\chi Pq$, $r_1^{(2)}$ solves $(2\partial+i\bar{k})r_1^{(2)}=(-2\partial)^{2}\chi P(qr_0)$. Since $||r_0||_{L^p(\Omega)}\sim \mathcal{O}(\frac{1}{|\bm{k}|})$, $||r_1^{(2)}||_{L^p(\Omega)}\sim \mathcal{O}(\frac{1}{|\bm{k}|})$, so we have
\begin{equation}
r_1=\frac{\chi Pq\chi Pq}{(i\bar{k})^2}+\mathcal{O}(\frac{1}{|\bm{k}|^3}).
\end{equation}
The rest $r_i$, $i\geq 2$ are solved up to the first order and $||r_i||_{L^p(\Omega)}\sim \mathcal{O}(\frac{1}{|\bm{k}|^{i+1}})$, so we get the expansion
\begin{equation}
r(\bm{x},\bm{k})=\frac{\chi Pq}{i\bar{k}}+\frac{(-2\partial+\chi Pq)\chi Pq}{(i\bar{k})^2}+\frac{b_3(\bm{x},\bm{k})}{(i\bar{k})^3}
\end{equation}

For $\tilde{\bm{\eta}}=\frac{1}{2}(-\bm{k}^\perp +i\bm{k})$, we can run the same process with some minor change and get the expansion expression for the corresponding correction term
\begin{equation}
\label{higherorder2}
\tilde{r}(\bm{x},\bm{k})=\sum_{j=1}^{n-1} \frac{\tilde{a}_j(\bm{x})}{(ik)^j}+\frac{\tilde{b}_n(\bm{x},\bm{k})}{(ik)^n},\quad
\tilde{a}_j=(-2\bar{\partial}+\chi \bar{P}q)^{j-1}\chi \bar{P}q
\end{equation} 

Now we can show Lemma \ref{smooth} using Lemma \ref{expand}. First we show the similar result for the Schr\"odinger equation, which is actually a generalization of the results in \cite{sugeneric}. To show that $\forall g\in span\{uv\}^\perp$, $g\chi_\Omega$ belongs to $H^s(\mathbf{R}^2)$ for $0<s<1$, they take CGO solutions
\begin{equation*}
u=e^{(\frac{1}{2}(\bm{k}^\perp +i\bm{k}))\cdot\bm{x}}(1+r(\bm{x},\bm{k})),\quad v=e^{(\frac{1}{2}(-\bm{k}^\perp +i\bm{k}))\cdot\bm{x}}(1+\tilde{r}(\bm{x},\bm{k})), 
\end{equation*}
for $|\bm{k}|$ sufficiently large. Plug these into the assumption $\int_\Omega guv\, d\bm{x}=0$, they get
\begin{equation}
\label{ghat}
\int_\Omega ge^{i\bm{k}\cdot\bm{x}}(1+r+\tilde{r}+r\tilde{r})\, d\bm{x}=0\quad
\Rightarrow\quad \widehat{g\chi_\Omega}(-\bm{k})=\int_\Omega ge^{i\bm{k}\cdot\bm{x}}(r+\tilde{r}+r\tilde{r})\, d\bm{x}.
\end{equation}
Then they use the expansion (\ref{higherorder}), (\ref{higherorder2}), up to the second order, to showed that $|\bm{-k}|^{s}\widehat{g\chi_\Omega}(-\bm{k})$ is in $L^2(|\bm{k}|>R)$ which implies $g\chi_\Omega\in H^s(\mathbf{R}^2)$ since only the large $\bm{k}$ behaviour of $\widehat{g\chi_\Omega}$ matters. They stopped here since it was enough for their purpose, while we can actually prove $g\chi_\Omega\in H^s$ for any fixed $s>0$ by expanding (\ref{higherorder}), (\ref{higherorder2}) to higher order plus a bootstrap process. 

For instance, if we expand $r$ and $\tilde{r}$ up to the third order, then (\ref{ghat}) becomes
\begin{equation*}
\begin{split}
&\widehat{g\chi_\Omega}(-\bm{k})=\frac{1}{i\bar{k}}\widehat{a_1 g\chi_\Omega}(-\bm{k})+\frac{1}{ik}\widehat{\tilde{a}_1 g\chi_\Omega}(-\bm{k})+\\
&\frac{1}{(i\bar{k})^2}\widehat{a_2 g\chi_\Omega}(-\bm{k})+\frac{1}{(ik)^2}\widehat{\tilde{a}_2 g\chi_\Omega}(-\bm{k})
+\frac{1}{(i\bar{k})(ik)}\widehat{a_1\tilde{a}_1 g\chi_\Omega}(-\bm{k})+\\
&\frac{1}{(i\bar{k})(ik)^2}\widehat{a_1\tilde{a}_2 g\chi_\Omega}(-\bm{k})+\frac{1}{(i\bar{k})^2(ik)}\widehat{a_2\tilde{a}_1 g\chi_\Omega}(-\bm{k})+\frac{1}{(i\bar{k})^2(ik)^2}\widehat{a_2\tilde{a}_2 g\chi_\Omega}(-\bm{k})+\\
&\frac{1}{(i\bar{k})^3}\int_\Omega g(\bm{x})b_3(\bm{x},\bm{k})(1+\tilde{r}(\bm{x},\bm{k}))e^{i\bm{k}\cdot\bm{x}}\, d\bm{x}
+\frac{1}{(ik)^3}\int_\Omega g(\bm{x})\tilde{b}_3(\bm{x},\bm{k})(1+r(\bm{x},\bm{k}))e^{i\bm{k}\cdot\bm{x}}\, d\bm{x},\\
\end{split}
\end{equation*}
for $|\bm{k}|$ large and from here we can show that for $0<s<2$, $|\bm{-k}|^{s}\widehat{g\chi_\Omega}(-\bm{k})\in L^2(|\bm{k}|>R)$, so $g\chi_\Omega\in H^s(\mathbf{R}^2)$. 
\begin{itemize}
\item The terms in the first line, like $\frac{|\bm{-k}|^s}{i\bar{k}}\widehat{a_1 g\chi_\Omega}(-\bm{k})$, belongs to $L^2(|\bm{k}|>R)$ is based on the fact that it have already been shown that $g\chi_\Omega\in H^t(\mathbf{R}^2)$ for $0<t<1$, $a_1$ is a regular enough function(depends on the regularity of the extension of $q$) with compact support, then $a_1g\chi_\Omega\in H^t(\mathbf{R}^2)$, which means $|\bm{-k}|^t\widehat{a_1g\chi_\Omega}(-\bm{k})\in L^2(|\bm{k}|>R)$.
\item The terms in the second and third line, like $|\bm{-k}|^{s}\frac{1}{(i\bar{k})^2}\widehat{a_2 g\chi_\Omega}(-\bm{k})$, belongs to $L^2(|\bm{k}|>R)$ since all the Fourier transform is in $L^2(\mathbf{R}^2)$ while the rest like $|\bm{-k}|^{s}\frac{1}{(i\bar{k})^2}$ is bounded.
\item The terms in the last line, like $\frac{|\bm{-k}|^{s}}{(i\bar{k})^3}\int_\Omega g(\bm{x})b_3(\bm{x},\bm{k})(1+\tilde{r}(\bm{x},\bm{k}))e^{i\bm{k}\cdot\bm{x}}\, d\bm{x}$, belongs to $L^2(|\bm{k}|>R)$ since $\frac{|\bm{-k}|^{s}}{(i\bar{k})^3}\in L^2(|\bm{k}|>R)$ given that $0<s<2$ and the rest $\int_\Omega g(\bm{x})b_3(\bm{x},\bm{k})(1+\tilde{r}(\bm{x},\bm{k}))e^{i\bm{k}\cdot\bm{x}}\, d\bm{x}$ is bounded in $\bm{k}$ by Cauchy-Schwarz inequality,
\begin{equation*}
|\int_\Omega g(\bm{x})b_3(\bm{x},\bm{k})(1+\tilde{r}(\bm{x},\bm{k}))e^{i\bm{k}\cdot\bm{x}}\, d\bm{x}|\leq
||g||_{L^2(\Omega)}||b_3(\bm{x},\bm{k})||_{L^4(\Omega)}||1+\tilde{r}(\bm{x},\bm{k})||_{L^4(\Omega)}.
\end{equation*}
The $L^4(\Omega)$ norm here is in $\bm{x}$ variable but is uniformly bounded in $\bm{k}$, which finishes the proof.
\end{itemize}
In general, for any fixed $s>0$, we can show $g\chi_\Omega\in H^s(\mathbf{R}^2)$ by a bootstrap process like above by expanding $r$ and $\tilde{r}$ up to higher order, which require more regularity on the potential and its extension. We remark here that this argument can be used to fix the little flaw in \cite{buk}, where it started with the assumption that $g\in C^{1}(\bar{\Omega})$. 

Now when we come back to the conductivity equation, just take the CGO solutions obtained by the Liouville transform, i.e.
\begin{equation*}
u=\gamma^{-\frac{1}{2}}e^{(\frac{1}{2}(\bm{k}^\perp +i\bm{k}))\cdot\bm{x}}(1+r(\bm{x},\bm{k})),\quad 
v=\gamma^{-\frac{1}{2}}e^{(\frac{1}{2}(-\bm{k}^\perp +i\bm{k}))\cdot\bm{x}}(1+\tilde{r}(\bm{x},\bm{k})),
\end{equation*}
plug into the assumption $\int_\Omega f\nabla u\cdot \gamma\nabla v\,dx=\frac{1}{2}\int_\Omega f\nabla (\gamma \nabla uv)\,dx=0$ and get
\begin{equation}
\begin{split}
&\widehat{f\chi_\Omega}(-\bm{k})=-\int_\Omega fRe^{i\bm{k}\cdot\bm{x}}\,d\bm{x}+\\
&\frac{2i}{|\bm{k}|^2}\bm{k}\cdot\int_\Omega f\nabla Re^{i\bm{k}\cdot\bm{x}}\,d\bm{x}
-\frac{i}{|\bm{k}|^2}\bm{k}\cdot\int_\Omega f(1+R)\nabla \log{\gamma} e^{i\bm{k}\cdot\bm{x}}\,d\bm{x}\\
&\frac{1}{|\bm{k}|^2}(\int_\Omega f(1+R)\Delta \log{\gamma} e^{i\bm{k}\cdot\bm{x}}\,d\bm{x}
+\int_\Omega f\nabla R\cdot\nabla \log{\gamma} e^{i\bm{k}\cdot\bm{x}}\,d\bm{x}
+\int_\Omega f\Delta R e^{i\bm{k}\cdot\bm{x}}\,d\bm{x}),\\
\end{split}
\end{equation}
where $R=R(\bm{x},\bm{k})=r(\bm{x},\bm{k})+\tilde{r}(\bm{x},\bm{k})+r(\bm{x},\bm{k})\tilde{r}(\bm{x},\bm{k})$. The leading order in $|\bm{k}|$ is actually the same as what we get in (\ref{ghat}), so we can run the same arguments to prove $f\chi_\Omega\in H^s(\mathbf{R}^2)$. Remember when proving the similar smoothness result for the Schr\"odinger equation, we need some regularity for the potential function, which is now related to the conductivity by the Liouville transform as $q=\gamma^{-\frac{1}{2}}\Delta \gamma^{\frac{1}{2}}$. As a result, here we need the conductivity to be regular enough to run the above arguments. 

Now we are ready to prove Lemma \ref{den}.
\begin{pf}
The lemma is equivalent with that for any $f\in L^2(\Omega)$, s.t. $f\perp span\{\nabla u\cdot \bm{\gamma}\nabla v\}$, where $u$, $v$ are arbitrary solutions of (\ref{ce}), $f$ is identically $0$.

As mentioned above, we may assume the conductivity $\bm{\gamma}$ is isotropic. By Lemma \ref{smooth}, if the conductivity is regular enough, we have $f\chi_\Omega\in H^s(\mathbf{R}^2)$, for any fixed  positive $s$. By the Soblev embedding, $f\chi_\Omega$ belongs to $C^m(\mathbf{R}^2)$, $\forall m\in\mathbf{N}_+$, which implies $f\in C^m(\mathbf{R}^2)$ and vanishing to $m$-th order at the boundary.

In particular, we prove up to $f\chi_\Omega\in C^3(\mathbf{R}^2)$ so both $f$ and $\partial_{\bm{\nu}}f$ vanish at the boundary, then by the assumption that $f\perp span\{\nabla u\cdot \gamma\nabla v\}$, we have
\begin{equation}
\label{ibp}
0=\int_\Omega f\nabla u\cdot \gamma\nabla v\, d\bm{x}=\frac{1}{2}\int_\Omega f\nabla\cdot (\gamma\nabla(uv))\, d\bm{x}=\frac{1}{2}\int_\Omega \nabla\cdot (\gamma\nabla(f))uv\, d\bm{x}.
\end{equation}
Use the Liouville transform, any solution $u$ of (\ref{ce}) can be written as $u=\gamma^{-\frac{1}{2}}u_s$, where $u_s$ is solution to the related Schr\"odinger equation. So we have $\int_\Omega \frac{1}{\gamma}\nabla\cdot (\gamma\nabla(f))u_s v_s\, d\bm{x}=0$. Since $f\in C^3(\mathbf{R}^2)$ and $\gamma$ is regular enough, $\frac{1}{\gamma}\nabla\cdot (\gamma\nabla(f))\in C^1(\bar{\Omega})$, then it have been shown in \cite{buk} that $\frac{1}{\gamma}\nabla\cdot (\gamma\nabla(f))=0$ by using special constructed solutions and the stationary phase method.

Finally by the unique continuation result for elliptic equations, $\frac{1}{\gamma}\nabla\cdot (\gamma\nabla(f))=0$ together with both $f$ and $\partial_{\bm{\nu}}f$ vanish at the boundary imply that $f=0$. 
\end{pf}

We remark here that the methods used above to construct the CGO solutions can also be used to construct the solutions Bukhgeim used in \cite{buk}, i.e. $\forall \phi(\bm{x})=u(\bm{x})+iv(\bm{x})$ analytic(or $\tilde{\phi}=-u+iv$ anti-analytic), we construct solutions for the Schr\"odinger equation in the form $u(\bm{x})=e^{i\tau\phi(\bm{x})}(1+r(\bm{x},\tau))$(or $v=e^{i\tau\tilde{\phi}}(1+\tilde{r})$). Since in this case, we also have $\Delta e^{i\tau\phi(\bm{x})}=0$(or $\Delta e^{i\tau\tilde{\phi}}=0$), and $\nabla \phi(\bm{x})=\nabla u(\bm{x})+i\nabla v(\bm{x})=(\nabla v)^{\perp}+i\nabla v$(or $\nabla \tilde{\phi}=-(\nabla v)^{\perp}+i\nabla v$) due to the Cauchy-Riemann equation, which is quite similar with the CGO solution where we have $\Delta e^{\bm{\eta}\cdot\bm{x}}=0$, $\nabla e^{\bm{\eta}\cdot\bm{x}}=\bm{\eta}=\frac{1}{2}(\bm{k}^{\perp}+i\bm{k})$. But we can not get the similar expansion structure as we have in Lemma \ref{expand}.

\section{Determination of the Heat Parameters}
\label{rh}
With all the preparation above, we can now determine the heat parameters in this section. 

Assume $\bm{F}$ is the  diffeomorphism obtained in Corollary \ref{rc2}, then according to the results in Section \ref{coc} and the uniqueness for the Dirichlet boundary value problem for (\ref{ce}), we have
\begin{equation}
\begin{split}
u_2(\bm{x},t)=&u_1(\bm{F}^{-1}(\bm{x}),t)\\
S_2(\bm{x},t)=&\nabla u_2(\bm{x},t)\cdot \bm{\gamma}_2(\bm{x})\nabla u_2(\bm{x},t)=\nabla u_1(\bm{F}^{-1}(\bm{x}),t)\cdot \bm{F}_*\bm{\gamma}_1(\bm{x})\nabla u_1(\bm{F}^{-1}(\bm{x}),t)\\
=&\frac{1}{|D\bm{F}(\bm{F}^{-1}(\bm{x}))|}S_1(\bm{F}^{-1}(\bm{x}),t).
\end{split}
\end{equation}
On the other hand, using a change of coordinates to $\bm{\gamma}_1$, $\kappa_1$, $\bm{A}_1$ as in Section \ref{coc},
\begin{equation}
\begin{split}
&\tilde{\bm{\gamma}}_1(\bm{x})=\bm{F}_*\bm{\bm{\gamma}_1}(\bm{x})\triangleq\frac{D\bm{F} \bm{A} D\bm{F}^T}{|D\bm{F}|}\circ \bm{F}^{-1}(\bm{x})=\bm{\gamma}_2(\bm{x})\triangleq \bm{\gamma}(\bm{x}),\\
&\tilde{\kappa}_1(\bm{x})\triangleq |D\bm{F}(\bm{F}^{-1}(\bm{x}))|\kappa_1(\bm{F}^{-1}(\bm{x})),\,
\tilde{A}_1(\bm{x})=\bm{F}_*\bm{A_1}(\bm{x})\triangleq \frac{D\bm{F} \bm{A} D\bm{F}^T}{|D\bm{F}|}\circ \bm{F}^{-1}(\bm{x}),
\end{split}
\end{equation}
we know that $\Sigma_{\tilde{\bm{\gamma}}_1,\tilde{\kappa}_1,\tilde{\bm{A}}_1}=\Sigma_{\bm{\gamma}_1,\kappa_1,\bm{A}_1}=\Sigma_{\bm{\gamma}_2,\kappa_2,\bm{A}_2}$, $\tilde{u}_1(\bm{x})=u_2(\bm{x})$, $\tilde{S}_1(\bm{x},t)=S_2(\bm{x},t)\triangleq S(\bm{x},t)$. As we mentioned above, we take special boundary data which is separated in $\bm{x}$ and $t$, i.e. $f(\bm{x},t)=h(\bm{x})g(t)$ then $S(\bm{x},t)=(\nabla u(\bm{x})\cdot \bm{\gamma}(\bm{x})\nabla u(\bm{x}))g(t)$, where $u(\bm{x})$ is arbitrary solution to (\ref{ce}). By the density arguments from Section \ref{dense}, the linearity and continuity of the heat equation, we know that the out-coming heat flow are the same for the two systems (one with coefficients $\tilde{\kappa}_1$, $\tilde{\bm{A}}_1$, the other with $\kappa_2$, $\bm{A}_2$) when we replace the source term $S(\bm{x},t)$ in (\ref{he}) with $\bar{S}(\bm{x},t)=w(\bm{x})g(t)$, where $w(\bm{x})$ is arbitrary function in $L^2(\Omega)$. From here, we can conclude that
\begin{equation}
\kappa_2=\tilde{\kappa}_1, \quad \bm{A}_2=\tilde{\bm{A}}_1,
\end{equation}
which has been proved in \cite{hybrid}. We are not going to repeat the proof here but the main idea is to take special input in the form $S(\bm{x},t)=w(\bm{x})\delta(t)$, then solve the equations using the eigenfunction expansions.

\section{Summary}
In summary, we generalize the higher dimensional uniqueness result for a new hybrid method proposed in \cite{hybrid} to the two dimensional case. The main difficulty is to show a density arguments, which is proved in Section \ref{dense}. Since the two dimensional anisotropic Calder\'on's inverse problem is well understood, in this paper we allow the conductivity to be anisotropic and so we can only expect to have uniqueness up to a boundary-fixing diffeomorphism. The method actually doesn't provide any interior information, so we should not expect any improvement in stability while on the other hand, we can recover three coefficients all together and the use of electric boundary sources may be easier than controlling the temperature at the boundary. Further work may include requiring less regularity of the parameters, numerical reconstruction, improve the model since the boundary voltage is varying in time, etc.

\section{Acknowledgement}
The author would like to thank the HKUST Jockey Club Institute for Advanced Study for their support during fall semester 2015. The author was also partially supported by NSF.

\bibliography{hybrid}

\end{document}